\documentclass{amsart}

\usepackage{graphicx}
\usepackage{amssymb}
\usepackage{hyperref}

\theoremstyle{plain}
\newtheorem{thm}[subsection]{Theorem}

\newtheorem{lem}[subsection]{Lemma}
\newtheorem{prop}[subsection]{Proposition} 
\newtheorem{cor}[subsection]{Corollary}
\newtheorem{prps}[subsection]{Properties}

\theoremstyle{remark}
\newtheorem{rem}[subsection]{Remark}
\theoremstyle{definition}
\newtheorem{dfn}[subsection]{Definition} 
\newtheorem{ejm}[subsection]{Example}

\newtheorem{cvn}[subsection]{Convention}

\newcommand{\bc}{\mathbb{C}}
\newcommand{\bn}{\mathbb{N}}
\newcommand{\bbf}{\mathbb{F}}
\newcommand{\bq}{\mathbb{Q}}
\newcommand{\bz}{\mathbb{Z}}
\newcommand{\br}{\mathbb{R}}
\newcommand{\bp}{\mathbb{P}}
\newcommand{\be}{\mathbb{E}}
\newcommand{\ba}{\mathbb{A}}
\newcommand{\bb}{\mathbb{B}}
\newcommand{\bd}{\mathbb{D}}
\newcommand{\bt}{\mathbb{T}}

\newcommand{\ca}{\mathcal{A}}
\newcommand{\cc}{\mathcal{C}}
\newcommand{\cv}{\mathcal{V}}
\newcommand{\one}{{\mathbf 1}}

\DeclareMathOperator{\sing}{Sing}
\DeclareMathOperator{\Hom}{Hom}
\DeclareMathOperator{\Tors}{Tors}
\DeclareMathOperator{\rank}{rk}

\numberwithin{equation}{section}

\newcommand\rightmap[1]{\smash{\mathop{\ \mapsto\ }\limits^{#1}}}
\newcommand\enet[1]{\renewcommand\theenumi{#1}
\renewcommand\labelenumi{\theenumi}}

\title[Fundamental groups and pencils]%
{On the connection between fundamental groups and pencils with multiple fibers}

\author[E.~Artal]
{Enrique~Artal~Bartolo}
\address{
Departamento de Matem\'aticas, IUMA,
Facultad de Ciencias,
Universidad de Zaragoza,
c/ Pedro Cerbuna, 12,
50009 Zaragoza, Spain.}
\email{artal@unizar.es}
\urladdr{http://riemann.unizar.es/geotop/WebGeoTo/Profes/eartal/}

\author[J.I.~Cogolludo]
{Jos\'e~Ignacio~Cogolludo-Agust{\'\i}n}
\address{Departamento de Matem\'aticas, IUMA,
Facultad de Ciencias,
Universidad de Zaragoza,
c/ Pedro Cerbuna, 12,
50009 Zaragoza, Spain.}
\email{jicogo@unizar.es}
\urladdr{http://riemann.unizar.es/geotop/WebGeoTo/Profes/jicogo/}

\thanks{Partially supported by MTM2007-67908-C02-01.}

\date{\today}

\begin{document}
\maketitle


\section*{Introduction}
The study of the topology of complex projective (or quasiprojective)
smooth varieties depends strongly on the knowledge of the topology of the complement
of hypersurfaces in a projective space. Considering a projection, any smooth projective variety
is a covering of a projective space of the same dimension ramified along a hypersurface. These coverings
are measured by (finite index subgroups of) the fundamental group of the complement of the hypersurface.
Using Lefschetz-Zariski theory, if we take a \emph{generic} plane section the fundamental group of the 
complement does not change. As a consequence, for fundamental group purposes, one can restrict their 
attention to the study of complements of curves in the projective plane, as stated in the foundational 
paper by O.~Zariski~\cite{zr:29}.

The richness of coverings for a space depends on its fundamental group. This is why we are mostly interested 
in curves $C\subset\bp^2$ whose $\pi_1(\bp^2\setminus C)$ is non-abelian. The first known example is 
probably the curve formed by three lines $C:=L_1\cup L_2\cup L_3$ intersecting at one point~$P$. There is 
an easy way to compute this fundamental group; the pencil of lines through~$P$ is parametrized by $\bp^1$; 
this pencil induces an epimorphism of $\bp^2\setminus C$ onto $\bp^1\setminus\{p_1,p_2,p_3\}$ (the punctures
corresponding to the three lines). Moreover, this map is a locally trivial fibration (with fiber isomorphic 
to $\bc$) and hence $\pi_1(\bp^2\setminus C)\cong\pi_1(\bp^1\setminus\{p_1,p_2,p_3\})$, which is a free group 
of two generators.

The first known examples of irreducible curves whose fundamental groups are known to be non-abelian appeared
in~\cite{zr:29}. The first one corresponds to a hexacuspidal sextic, with its six cusps on a  conic; the 
equation of such a curve is of the form $f_2^3+f_3^2=0$, where $f_j$ is a homogeneous polynomial of degree~$j$. 
Its fundamental group is $\bz/2*\bz/3$; in \S\ref{sec-ex} we will see the relation between this group and the 
pencil generated by $f_2^3=0$ and $f_3^2=0$. This kind of examples have been generalized by various authors 
replacing $(2,3)$ by $(p,q)$. In the same paper, Zariski found the irreducible curve with smallest possible 
degree having a non-abelian fundamental group: the tricuspidal quartic. This example and many others appearing 
in the literature are also connected with pencils.

The precise connection with pencils can be stated as follows: a pencil defines a dominant morphism to a
quasi-projective curve, inducing an epimorphism at the level of fundamental groups. The multiplicities of the 
fibers of the pencil induce an orbifold structure on the quasi-projective group, and the map defines an 
epimorphism onto the orbifold fundamental group. When such an orbifold fundamental group is non-abelian,
then the original fundamental group has a surjection onto a non-abelian group. Such surjections coming 
from dominant maps will be referred to as \emph{geometric} surjections.

The tricuspidal quartic is the only irreducible curve of degree~$4$ with a non-abelian fundamental group.
The degree-five case was studied by A.~Degtyarev~\cite{deg:90}; he found exactly two irreducible quintics 
with non-abelian fundamental groups. One of them, also studied by the first author~\cite{ea:ta}, has an 
infinite fundamental group. In~\S\ref{sec-ex}, we will study its relationship with a pencil. The question
whether or not all non-abelian fundamental groups have a \emph{geometric} surjection onto an orbifold 
group naturally arises. All the examples studied, up to now, supported an affirmative answer to this question.

In this paper we will show an explicit example of a non-abelian fundamental group whose complement admits 
no geometric surjections. This curve is one of the quintics referred to in the previous paragraph, which 
will be called the \emph{projective Degtyarev curve} throughout this text.
As a brief description, the projective Degtyarev curve has exactly three singular points of type $\ba_4$;
its fundamental group is finite and non-abelian. In Proposition~\ref{prop-afin}, we prove that this group 
admits no geometric surjections. Once the group is computed, the proof is rather straightforward; it 
depends on the orders of the group and its abelianization and on the properties of orbifold groups. 

If we add a tangent line to one of the singular points of the projective Degtyarev curve, the 
complement of the union in $\bp^2$ is the complement of an affine curve, which will be called the 
\emph{affine Degtyarev curve}. This affine curve has an infinite non-abelian fundamental 
group and non-trivial characteristic varieties (see \S\ref{sec-orb} for the definition). Extending results of Arapura and others, it is known that irreducible components of positive dimension (for the fundamental group of a quasiprojective variety) are obtained as pull-back of irreducible components of characteristic varieties of orbifolds. A natural question arises: Is it also true for irreducible components of dimension~$0$ (isolated points)? Plenty of computations supported a positive answer: most quasiprojective groups satisfy the property 
for irreducible components of any dimension (see~\S\ref{sec-ex} for examples). The main Theorem~\ref{thm-main} of this paper shows that the fundamental group of the complement of the affine Degtyarev curve does not satisfy 
this property. This is the only known example, up to now.

The paper is organized as follows. In~\S\ref{sec-orb}, the concepts of orbifold and characteristic varieties 
are recalled, also some orbifold groups are studied. In~\S\ref{sec-ex}, we relate non-abelian fundamental groups of the complements of curves (which are known in the literature) with orbifold morphisms (via pencils of curves). In~\S\ref{sec-dgt}, we describe Degtyarev curves and, in order to obtain a prresentation for their 
fundamental groups, we compute a special braid monodromy. The fundamental groups are obtained 
in~\S\ref{sec-groups}, where also the main results of the paper are stated and proved. Finally, 
further properties of the affine Degtyarev curve are sketched in~\S\ref{sec-ell}.

\section{Orbifold groups and characteristic varieties}\label{sec-orb}

The fundamental groups of oriented Riemann surfaces have been extensively studied.
The fundamental group of a compact Riemann surface of genus~$g$ is
\[
\pi_g:=\left\langle a_i, b_i, 1\le i\le g\left| \prod_{i=1}^{g} a_ib_ia_i^{-1}b_i^{-1}\right.\right\rangle.
\]
If $C$ is a surface with genus~$g$ and $k>0$~punctures then its fundamental group
is free of rank~$2 g+k-1$. We are going to extend this family by considering \emph{orbifold groups}.

\begin{dfn}
An \emph{orbifold} $X_\varphi$ is a quasiprojective Riemann surface~$X$
with a function $\varphi:X\to\bn$ with value~$1$ outside a finite
number of points.
\end{dfn}

We may think that a neighborhood of a point $P\in X_\varphi$ such that
$\varphi(P)=n$ is the quotient of a disk (centered at~$P$) by a rotation of
angle $\frac{2\pi}{n}$. We will consider that a loop around $P$ is trivial if its lifting bounds a disk.
Following this idea, we define orbifold fundamental groups.

\begin{dfn}\label{dfn-group-orb}
For an orbifold $X_\varphi$, let $p_1,\dots,p_n$ the points such that
$m_j:=\varphi(p_j)>1$. Then, the \emph{orbifold fundamental group} of $X_\varphi$ is
$$
\pi_1^{\text{\rm orb}}(X_\varphi):=\pi_1(X\setminus\{p_1,\dots,p_n\})/\langle\mu_j^{m_j}=1\rangle
$$
where $\mu_j$ is a meridian of $p_j$. We denote $X_\varphi$ by $X_{m_1,\dots,m_n}$.
\end{dfn}

\begin{ejm}
If $X$ is a compact surface of genus~$g$ and type $X_{m_1,\dots,m_n}$, then
\[
\pi_1^{\text{\rm orb}}(X_\varphi)=\left\langle
a_1,\dots,a_g,b_1,\dots,b_g,\mu_1,\dots,\mu_n\left\vert
\prod_{i=1}^g [a_i,b_i]=\prod_{j=1}^n\mu_j,\ \underset{j=1,\dots,10}{\mu_j^{m_j}=1}\right.\right\rangle.
\]
If $X$ is not compact and $\pi_1(X)$ is free of rank~$r$, then
\[
\pi_1^{\text{\rm orb}}(X_\varphi)=\left\langle
a_1,\dots,a_r,\mu_1,\dots,\mu_n\left\vert
\underset{j=1,\dots,10}{\mu_j^{m_j}=1}\right.\right\rangle.
\]
\end{ejm}

\begin{dfn}
A dominant algebraic morphism $\rho:Y\to X$ defines an \emph{orbifold morphism} $Y\to X_\varphi$
if for all $p\in X$, the divisor $\rho^*(p)$ is a $\varphi(p)$-multiple.
\end{dfn}

\begin{prop}\label{prop-orb}
Let $\rho:Y\to X$ define an \emph{orbifold morphism} $Y\to X_\varphi$. Then $\rho$ induces a morphism
$\rho_*:\pi_1(Y)\to\pi_1^{\text{\rm orb}}(X_\varphi)$. Moreover, if the generic fiber is connected, then
$\rho_*$ is surjective.
\end{prop}

\begin{proof}
Let $M_\varphi:=\{x\in X\mid \varphi(x)>1\}$; we consider the restriction mapping $\tilde{\rho}:=\rho_|:
Y\setminus\rho^{-1}(M_\varphi)\to X\setminus M_\varphi$. This map induces a morphism
$\tilde{\rho}_*:\pi_1(Y\setminus\rho^{-1}(M_\varphi))\to\pi_1(X\setminus M_\varphi)$ fitting in the following commutative diagram:
\begin{equation*}
\begin{matrix}
\pi_1(Y\setminus\rho^{-1}(M_\varphi))&\overset{\tilde{\rho}_*}{\longrightarrow}&\pi_1(X\setminus M_\varphi)\\
i_*\downarrow&&\downarrow j_*\\
\pi_1(Y)&\overset{\rho_*}{\longrightarrow}&\pi_1(X).
\end{matrix}
\end{equation*}
The vertical mappings are induced by the inclusions. They are both surjective; the kernel of $j_*$ is generated
by the meridians of the points in $M_\varphi$ while the kernel of $i_*$ is generated by the meridians of the
irreducible components of $\rho^{-1}(M_\varphi)$, i.e., the components of the pull-back divisor $\rho^*(M_\varphi)$.

Let us consider an irreducible component $D$ of $\rho^*(M_\varphi)$ such that $\rho(D)=:x\in M_\varphi$. Let
$n:=\varphi(x)$; note that the multiplicity $m_D$ of $D$ in  $\rho^*(M_\varphi)$ is a multiple of~$n$.
We can interpret $m_D$ as follows. If $\mu_D$ denotes a meridian of $D$, then there is a meridian $\mu_x$ of~$x$
such that $\tilde{\rho}_*(\mu_D)=(\mu_x)^{m_D}$. Following Definition~\ref{dfn-group-orb}, it is easily seen
that $\tilde{\rho}_*$ factorizes through a morphism (also called $\rho_*$) 
$\pi_1(Y)\to\pi_1^{\text{\rm orb}}(X_\varphi)$.

The above argument also works if one replaces $M_\varphi$ by a finite set $M\supseteq M_\varphi$. 
In particular, one can choose $M$ to be the bifurcation locus of $\rho$, i.e., the mapping is a differentiable locally trivial fibration outside $M$. If the fiber is generically connected, the long exact homotopy sequence 
of this fibration implies the surjectivity of $\tilde{\rho}_*$ (for $M$). The result follows.
\end{proof}

\begin{dfn}
A fundamental group $G:=\pi_1(Y)$ is said to \emph{posses a geometric surjection} if $Y$ possesses an 
orbifold morphism $Y\to X_\varphi$ whose generic fiber is connected, and such that 
$\pi_1^{\text{\rm orb}}(X_\varphi)$ is non-abelian.
\end{dfn}

We recall the notion of characteristic varieties and its relationship with orbifolds.
We focus our attention on the characteristic varieties of quasiprojective manifolds, though they
can be defined in general and depend only on the fundamental group. Let $X$ be a connected 
topological space $X$, having the homotopy type of a finite $CW$-complex, and let $G:=\pi_1(X,x_0)$,
$x_0\in X$ which will be omitted if it is not necessary. Recall that the space of characters of $G$ is
\begin{equation}
\label{eq-torus}
H^1(X;\bc^*)=\Hom(H_1(X;\bz),\bc^*)=\Hom(\pi_1(X),\bc^*)=:\bt_G.
\end{equation}
\begin{rem}
Since $G$ is finitely generated, then it is also the case for $H_1(X; \bz)$.
Let $n:=\rank H_1(X; \bz)$ and $\Tors_G$ be the torsion subgroup of $H_1(X; \bz)$. Then
$\bt_G$ is an abelian complex Lie group with $|\Tors_G|$ connected components (each one is isomorphic
to $(\bc^*)^n$) satisfying the following exact sequence:
\[
1\to \bt_G^\one\to\bt_G\to \Tors_G\to 1,
\]
where $\bt_G^\one$ is the connected component containing the trivial character $\one$.
\end{rem}

For a character $\xi\in\bt_G$, we can construct a local system of coefficients $\bc_\xi$ over $X$.

\begin{dfn}
\label{def-char-var}
The $k$-th~\emph{characteristic variety}~of $X$ is the subvariety 
of $\bt_G$, defined by:
\[
\cv_{k}(X)=\{ \xi \in \bt_G\: |\:\dim H^1(X,\bc_{\xi}) \ge k \}, 
\]
where $H^1(X,\bc_{\xi})$ is the cohomology with coefficients 
in the local system $\xi$. In some cases we will use the notation $\cv_k(G)$.
\end{dfn}

The following result is straightforward.

\begin{prop}
Let $\varphi:G\to H$ be a group epimorphism. Then $\varphi^*$ induces injections $\bt_H\equiv\varphi^*\bt_H\hookrightarrow\bt_G$
and $\cv_j(H)\equiv\varphi^*\cv_j(H)\hookrightarrow\cv_j(G)$.
\end{prop}

\begin{rem}\label{complex}
Let us explain how to compute these invariants. For the sake of simplicity, the twisted homology, instead 
of the cohomology, will be computed. Let us consider a finite $CW$-complex homotopy equivalent to $X$. Let 
$\pi:\tilde{X}\to X$ be the maximal abelian covering. Note that $\tilde{X}$ inherits a $CW$-complex structure. 
The group of automorphisms of $\pi$ is $H_1(X;\bz)$. The action of this Abelian group endows the chain complex
$C_*(\tilde{X};\bc)$ with a module structure over the ring $\Lambda:=\bz[H_1(X;\bz)]$. The differentials of the 
complex are $\Lambda$-homomorphisms. Moreover, $C_*(\tilde{X};\bc)$ is a free $\Lambda$-module of finite rank.
If we fix a character $\xi$, $\bc$ has a natural $\Lambda$-module structure which is denoted by $\bc_\xi$ 
(as the local system of coefficients). The twisted homology of $X$ is the homology of the 
$C_*(X;\bc)^\xi:=C_*(\tilde{X};\bc)\otimes_\Lambda\bc_\xi$. Following this interpretation, it is not difficult 
to prove that the characteristic varieties are algebraic subvarieties of $\bt_G$, defined with integer equations.

This $i$-th jumping loci of $C_*(\tilde{X};\bc)$ with respect to $\underline{\ \ }\, \otimes_\Lambda\bc_\xi$ 
can also be viewed as the zero locus of the $i$-th Fitting ideal of $H_1(\tilde X;\bc)$ or, analogously, 
the support of the module $\wedge^i H_1(\tilde X;\bc)$ over the ring $\Lambda$ (see~\cite{li:01}).
\end{rem}

Following the theory developed by 
various authors (Beauville~\cite{Be}, Arapura~\cite{ara:97}, Simpson~\cite{Si1}, Budur~\cite{Bu}, 
Delzant~\cite{Delzant}, Dimca~\cite{Di4}), the structure of characteristic varieties for 
quasiprojective manifolds) can be stated as follows. 

\begin{thm}[\cite{ACM-prep}]\label{thm-orb}
Let $\Sigma$ be an irreducible component of $\cv_k(G)$, $k\geq 1$. Then one of the two following
statements holds:
\begin{itemize}
\item There exists a surjective orbifold morphism $\rho:X\to C_\varphi$
and an irreducible component $\Sigma_1$ of $\cv_k(\pi_1^{\rm\text{orb}}(C_\varphi))$
such that $\Sigma=\rho^*(\Sigma_1)$.
\item $\Sigma$ is an isolated torsion point.
\end{itemize}
\end{thm}

\begin{rem}
In general, both $G$ and its characteristic varieties are difficult to compute. For the complement of 
hypersurfaces in a projective space, Libgober~\cite{li:01} gave an alternative way of computing most 
components of the characteristic varieties from algebraic properties of the hypersurface without computing~$G$.
\end{rem}

\begin{rem}
Characteristic varieties can also be understood from Alexander-invariant point of view. Following Theorem~\ref{thm-orb}, characteristic varieties are determined by finite-index abelian coverings.
\end{rem}

We compute the invariants for some orbifold groups.

\begin{prop}\label{prop-2510}
Let $G$ be the orbifold group of~$\bp^1_{2,5,10}$.
Then $G$ is a semidirect product of the fundamental group of a compact surface of genus~$2$ and $\bz/{10}$.
The torus $\bt_G$ is $\mu_{10}$, the group of $10$-th roots of unity, $\cv_1(G)$ consists of the primitive $10$-th roots of unity and $\cv_2(G)=\emptyset$.
\end{prop}

\begin{proof}
Let us consider the short exact sequence associated with the abelianization map ($\bz/10:=\langle t\mid t^{10}=1\rangle$ is $G/G'$).
This sequence corresponds to a uniformization covering of the orbifold, and using Riemann-Hurwitz one obtains a compact Riemann surface of genus $2$. Since the exact sequence splits, we have a semidirect-group structure.

In order to compute $\cv_1(G)$ we follow the construction outlined in Remark~\ref{complex}, 
applied to the $CW$-complex associated with the presentation of $G$ given by
$\langle x,y\mid x^2=y^5=(x y)^{10}=1\rangle$.
Let us denote $p$ the unique $0$-cell, $x,y$ the $1$-cells and $A,B,C$ the $2$-cells (corresponding to the 
relations in the given order). Let us fix a character $\xi\in \bt_G$. It is clear that $\one\notin\cv_1(G)$. 
We can assume that $\zeta:=\xi(t)\neq 1$. The complex $C_*(X;\bc)^\xi$ is given by 
$$
0\longrightarrow\bc^3\overset{\partial_2}{\longrightarrow}\bc^2\overset{\partial_1}{\longrightarrow}\bc\longrightarrow 0.
$$
The matrix for $\partial_1$ is
$
\left(
\begin{smallmatrix}
\zeta^5-1\\
\zeta^2-1
\end{smallmatrix}
\right)$.
In particular, $\dim\ker\partial_1=1$. The matrix for $\partial_2$ equals
$$
\begin{pmatrix}
\zeta^5+1& 0&\dfrac{\zeta^{10}-1}{\zeta-1}\\
0&\zeta^8-\zeta^6+\zeta^4-\zeta^2+1&\zeta^5\dfrac{\zeta^{10}-1}{\zeta-1}
\end{pmatrix}
$$
In order to have non-trivial homology, this matrix must vanish and this happens only when $\zeta$ is a primitive $10$-th root of unity.
\end{proof}

\begin{prop}\label{prop-2255}
Let $G$ be the orbifold group of $\bp^1_{2,2,5,5}$.
Then $G$ is an extension of $\bz/{10}$ by the fundamental group of a compact surface of genus~$4$.
The torus $\bt_G$ is $\mu_{10}$, the group of $10$-th roots of unity, and both~$\cv_1(G)$ and $\cv_2(G)$ consist 
of the primitive $10$-th roots of unity.
\end{prop}

\begin{proof}
The short exact sequence associated with the abelianization map ($G/G'=\bz/10$) corresponds to 
a uniformization covering of the orbifold, and using Riemann-Hurwitz one obtains a compact Riemann 
surface of genus~$4$. 

We compute the characteristic varieties as in the proof of Proposition~\ref{prop-2510} for the presentation of $G$ given by $\langle x,y,z\mid x^5=y^5=z^2=(x y z)^{2}=1\rangle$.
Let us denote $p$ the unique $0$-cell, $x,y,z$ the $1$-cells and $A,B,C,D$ the $2$-cells (corresponding to the relations in the given order). Let us fix a character $\xi \in \bt_G$. It is clear that $\one\notin\cv_1(G)$. We can assume that
$\zeta:=\xi(t)\neq 1$. The complex $C_*(X;\bc)^\xi$ is given by 
$$
0\longrightarrow\bc^4\overset{\partial_2}{\longrightarrow}\bc^3\overset{\partial_1}{\longrightarrow}\bc\longrightarrow 0.
$$
The matrix for $\partial_1$ is the transposed of
$
\left(
\begin{smallmatrix}
\zeta^2-1&
\zeta^2-1&
\zeta^5-1
\end{smallmatrix}
\right)$.
In particular, $\dim\ker\partial_1=2$. The matrix for $\partial_2$ equals
$$
\begin{pmatrix}
\zeta^8-\zeta^6+\zeta^4-\zeta^2+1&0& 0&\zeta^5+1\\
0&\bar\zeta^8-\bar\zeta^6+\bar\zeta^4-\bar\zeta^2+1&0&\zeta(\zeta^5+1)\\
0&0&\zeta^5+1&\zeta^5+1
\end{pmatrix}
$$
In order to have non-trivial homology, this matrix must have rank less than~$2$ and this happens 
only when $\zeta$ is a primitive $10$-th root of unity. Moreover, in that case, the matrix vanishes.
\end{proof}

\section{Examples}\label{sec-ex}

In this section, we will present a collection of examples of curves with non-abelian fundamental groups
and geometric surjections and its relationship with characteristic varieties. 


\begin{rem}
\label{rem-lines}
If $Y:=\bp^2\setminus \cc$ admits an orbifold morphism $Y\to X_\varphi$, then the non-singular compactification 
$\bar X$ of $X$ is $\bp^1$.
\end{rem}

\begin{rem}\label{rem-lines2}
The easiest examples of curves with non-abelian fundamental groups and geometric surjections come from 
hyperplane (or line) arrangements. If a line arrangement $\ca$ has a point $P$ of multiplicity~$k\geq 3$, 
then the pencil of lines through~$P$ defines a morphism $\rho:\bp^2\setminus\ca\to X$, where $X$ is a 
$k$-punctured projective line. We have an epimorphism $\rho_*:\pi_1(\bp^2\setminus\ca)\to\pi_1(X)$ and the 
latter is a free group of rank~$k-1$ (hence non abelian). 
\end{rem}

The following result is well known for specialists.

\begin{prop}
The following three assertions are equivalent:
\begin{enumerate}
\item\label{prop-nabel} The group $\pi_1(\bp^2\setminus\ca)$ is non abelian,
\item\label{prop-m3} The arrangement $\ca$ has a point of multiplicity at least~$3$,
\item\label{prop-geomsurj} The group $\pi_1(\bp^2\setminus\ca)$ has a geometric surjection.
\end{enumerate}
\end{prop}

\begin{proof}
By the remark above, it is obvious that~(\ref{prop-m3}) implies~(\ref{prop-nabel}) and (\ref{prop-geomsurj}).
Also, by definition, (\ref{prop-geomsurj}) implies~(\ref{prop-nabel}). Hence it is enough to 
prove that~(\ref{prop-nabel}) implies~(\ref{prop-m3}). Note that, if~(\ref{prop-m3}) does not hold, then $\ca$ 
is an arrangement in general position. Either we choose a particular example (e.g. a real arrangement) and a 
braid monodromy argument implies immediately the abelianity or we use Hattori's topological description of 
arrangements of hyperplanes in general position.
\end{proof}

The argument used in Remark~\ref{rem-lines2} can be easily generalized when, instead of
considering three (or more) incident lines, one considers three (or more) fibers of any pencil of curves in 
$\bp^2$. Of course, any such example corresponds to curves with at least three irreducible components. 
The notion of orbifold allows for wider generalizations of this concept to curves with any number of 
irreducible components (for example to irreducible curves).

As it was stated in the Introduction, the first example of this kind is rather old, see~\cite{zr:29}. 
Let us consider a conic $C_2$ of equation $f_2=0$ and a cubic $C_3$ of equation $f_3=0$; we assume they do not 
have common components and they are not \emph{multiple} lines. Let $C$ be a curve of equation $f_2^3-f_3^2$. 
Note that the mapping $\rho:\bp^2\to \bp^1\setminus\{[1:1]\}$ 
given by $[x:y:z]\mapsto [f_2(x,y,z)^3:f_3(x,y,z)^2]$ is well defined 
(all the base points of the pencil are in $C$) and surjective. This mapping induces an orbifold map onto a $1$-punctured Riemann sphere with two \emph{orbifold} points of multiplicities $2$ and $3$ 
(at $[0:1]$ and $[1:0]$ respectively). Thus according to Proposition~\ref{prop-orb}, one obtains an epimorphism
$\pi_1(\bp^2\setminus C)$ onto $\bz/2*\bz/3$.

\begin{prop}
Let $G$ be the orbifold fundamental group of $\bc_{2,3}$. 
Then, $T_G=\mu_6$, $\cv_1(G)$ consists of the $6$-th primitive roots of unity and $\cv_2(G)=\emptyset$. 
In particular, any curve with equation $f_2^3-f_3^2=0$ has non-trivial characteristic varieties.
\end{prop}

\begin{rem}
For generic choices of $f_2$ and $f_3$ this epimorphism is in fact an isomorphism (this is actually the 
case originally considered by Zariski in~\cite{zr:29}). However, this is not the case, for instance, when 
$C$ is reducible (since $b_1(\bp^2\setminus C)>1$). Even if $C$ is irreducible one may also have not an 
isomorphism for several reasons: either there are few non-generic fibers in the pencil (e.g., a sextic 
with six cusps and four ordinary nodes) or there are several pencils (a sextic with nine cusps).
\end{rem}

These examples can be generalized if we replace $(2,3)$ by any coprimes $(p,q)$, see Oka~\cite{oka:75},
N\'emethi~\cite{nem:87} and Dimca~\cite{dim:li}. In such cases, the fundamental group of a generic curve with equation $f_p^q+f_q^p=0$ is $\bz/p*\bz/q$. Also Zariski~\cite{zr:29} considered another interesting example where 
the target orbifold is compact.

Let us consider the tricuspidal quartic $C_4$ with equation $f_4=0$. It is not hard to prove that
we can choose 
\begin{equation}\label{eq-tricusp}
f_4:=x^2 y^2+ y^2 z^2+ x^2 z^2- 2 x y z (x+y+z).  
\end{equation}
and $\sing(C_4)=\{[1:0:0],[0:1:0],[0:0:1]\}$.
The curve $C_4$ is parametrized by 
\begin{equation}\label{par-tricusp}
[t:s]\mapsto[t^2s^2: (t-s)^2s^2: t^2 (t-s)^2];
\end{equation}
and its singular points correspond to 
$[t:s]=[0:1],[1:1],$ and $[1:0]$. Let $P\in C_4$ be a smooth point with parameter~$\alpha\equiv[\alpha:1]$ and 
let $L$ be the tangent line to $C_4$ at $P$, with equation $f_1=0$, where 
\[
f_1:=(\alpha-1)^3 x-\alpha^3 y+z. 
\]
Let $C_2$ be the conic passing through 
the singular points of $C_4$ and tangent to $C_4$ at $P$; since we have imposed five (non-degenerate) conditions, 
such a conic is unique. As before, let $f_2=0$ be the equation of $C_2$, where
\[
f_2:=\alpha(\alpha-1) x y-(\alpha-1) x  z+\alpha y z.
\]
We consider now a cubic $C_3$ having 
a nodal point at~$P$ (one of the branches tangent to $C_4$ at~$P$) and tangent to $C_4$ at the three cuspidal 
points. Counting the conditions it is easy to prove that only one such cubic exists, with equation $f_3=0$, where
\[
f_3:=- \left( \alpha-2 \right)  \left( 2\,\alpha-1 \right)  \left( \alpha+1
 \right) xyz-{\alpha}^{3}x{y}^{2}-x{z}^{2}- \left( \alpha-1 \right) ^{
3}{x}^{2}y+y{z}^{2}+ \left( \alpha-1 \right) ^{3}{x}^{2}z+{\alpha}^{3}
{y}^{2}z.
\]

\begin{lem}\label{lem-tricusp}
$f_4 f_1^2 =f_3^2-4 f_2^3$.
\end{lem}

\begin{prop}
\label{prop-pencil-3cusp}
The fundamental group of $\bp^2\setminus C_4$ possesses a geometric surjection onto $\bp^1_{2,2,3}$.
\end{prop}

\begin{rem}
Zariski proved in~\cite{zr:29} that $\pi_1(\bp^2\setminus C_4)$ is non-abelian of order~$12$. The above mapping 
induces a central extension of $\bd_6$ (dihedral group of order~$6$) whose kernel is cyclic of order~$2$. Note 
that we also have an epimorphism from $\pi_1(\bp^2\setminus (C_4\cup L))$ onto the orbifold group of a 
$1$-punctured Riemann sphere with two multiple points $(2,3)$; for a generic $P$ it is possible to prove that 
$\pi_1(\bp^2\setminus (C_4\cup L))$ equals $\bb_3$. There is a particular case corresponding to the bitangent line; 
in this case there are two such mappings and  $\pi_1(\bp^2\setminus (C_4\cup L))$ is the Tits-Artin group of a
triangle.
\end{rem}

In~\cite{deg:90}, Degtyarev proved that only two irreducible curves of degree~$5$ have non-abelian fundamental group. One of them is extensively studied in~\S\ref{sec-dgt}. The other one was also studied by the first author in~\cite{ea:ta}. It is a rigid curve with one point of type~$\ba_6$ and three cuspidal points (it is the dual curve of the quartic with one~$\ba_6$).
Let $C_5$ be this curve (with equation~$f_5=0$). The way to prove that this group is non-abelian in~\cite{ea:ta} is to show that there is an epimorphism from an actual presentation of $\pi_1(\bp^2\setminus C_5)$ onto the triangle group of type~$2,3,7$; this is the orbifold
group of $\bp^1_\varphi$ with three multiple points of these orders. 
In fact, one has the following:

\begin{prop}
The fundamental group of $\bp^2\setminus C_5$ possesses a geometric surjection onto $\bp^1_{2,3,7}$
\end{prop}

\begin{proof}
We start with a pencil as in Lemma~\ref{lem-tricusp}, where $L$ is the bitangent. 
We are going to consider the Cremona transformation $\rho:\bp^2\dashrightarrow\bp^2$ associated with 
the net of conics having three infinitely near points in common with $C$ at $P$, the singular point of 
type $\ba_6$. Let us describe this transformation. We blow-up this three infinitely near points and we 
obtain a rational surface~$X$ with a morphism $\sigma_1:X\to\bp^2$. Let us denote the three exceptional 
components (in order of appearance) by $E_1$, $E_2$, and $T$, and finally the tangent line of $C$ at $P$
by $L$, see Figure~\ref{fig-cremona1}. 

\begin{cvn}\label{cvn-str}
For birational morphisms, we keep the notation of a curve for its strict transform unless otherwise stated.
\end{cvn}
\begin{figure}
\centering
\includegraphics[scale=.5]{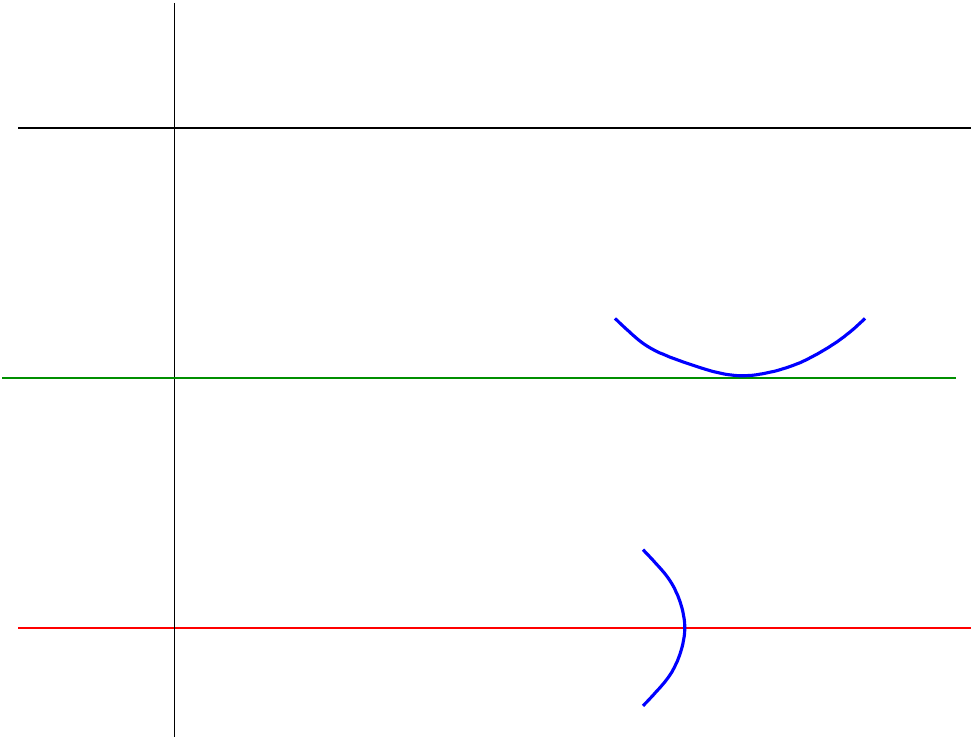}
\begin{picture}(0,0)
\put(-54,30){$C$}
\put(-54,60){$C$}
\put(-155,15){$L$}
\put(-155,90){$E_1$}
\put(-132,65){$E_2$}
\put(-90,55){$T$}
\end{picture}
\caption{Cremona transformation}
\label{fig-cremona1}
\end{figure}
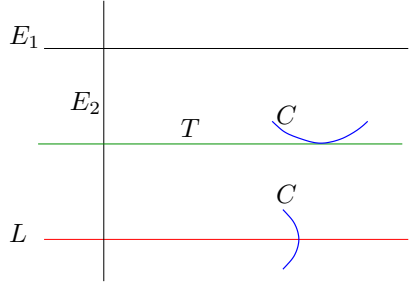
In $X$ we have $E_1\cdot E_1=-2$ and $E_1\cdot E_1=T\cdot T=L\cdot L=-1$. Since we can interchange $L$ and $T$,
we consider the birational morphism $\sigma_2:X\to\bp^2$ obtained as the composition of the contractions of $L$, $E_2$ and $E_1$; the contracted surface is rational with Euler characteristic~$3$ and hence it is a projective plane.
It is not hard to prove that $\rho=\sigma_2^{-1}\circ\sigma_1$. Let us denote $\tilde{C}:=\rho(C)$. Note that
$\tilde{C}$ is a tricuspidal quartic and $T$ is its unique bitangent line, one point $\tilde{P}$ comes from  the infinitely near point of $C$ at $P$ and the other one $Q$ comes from the other intersection point of $C$ and $L$.

We consider the pencil defined by the orbifold map of Proposition~\ref{prop-pencil-3cusp}, 
where the base point is~$\tilde{P}$. Let $C_3$ be the cubic such that $2 C_3$ is in the pencil. 
Following $C_3$ by $\sigma_2$ and $\sigma_1$, 
$C_6:=\rho^*(C_3)$ is a sextic with only one singular point at $P$ (with two branches, one of type $\ba_6$ and a smooth branch with maximal contact with the singular branch). With the same ideas, if $C_2$ is the conic such that $2 C_3$ is in the pencil, then $C_4:=\rho^*(C_2)$ is a quartic with an $\ba_6$ singular point at $P$.

Finally $\rho^*(\tilde{C}+2T)=C+7 L$. We have a pencil of degree~$12$ containing the fibers $2 C_6$, $3 C_4$ and $C+7 L$. This pencil produces the desired morphism.
\end{proof}

One can find more examples in the literature: Degtyarev~\cite{deg:90}, Flenner-Za{\u\i}denberg~\cite{fz}, and Tono~\cite{tono:01}. In what follows, the last two families will be described. We start with some definitions.

\begin{dfn}
A \emph{Hirzebruch} surface is a rational surface~$X$ with a morphism $\pi:X\to\bp^1$ which is a holomorphic (or algebraic) fibration with fiber~$\bp^1$. Such a surface is either $\Sigma_0:=\bp^1\times\bp^1$ or it has a unique section~$S_n$ with negative self-intersection $-n$, $n>0$; in that case $\pi$ is unique and $X$ is denoted by
$\Sigma_n$ ($\Sigma_1$ is the blowing-up of one point in $\bp^2$).
\end{dfn}

For any Hirzebruch surface~$X$ there is a family of birational maps which are called 
\emph{elementary Nagata transformations}. They are obtained as follows. Let us consider 
$\pi:X\to\bp^1$, $P\in X$ and $F:=\pi^{-1}(\pi(P))$; we consider the blowing-up 
$\sigma:\hat{X}\to X$ of $P$, with exceptional component~$\tilde{F}$. Since $(F\cdot F)_X=0$, 
we have that $(F\cdot F)_{\hat{X}}=-1$. By Castelnuovo criterion, we can blow down $F$
and we obtain a new Hirzebruch surface $\tilde{X}$ where $\tilde{F}$ is a fiber. 

\begin{dfn}
An elementary Nagata transformation is said to be \emph{positive} (resp. \emph{negative}) if $P$ 
belongs (resp. does not belong) to a section with non-positive self-intersection. For a positive one, 
one goes from $\Sigma_n$ to $\Sigma_{n+1}$; for a negative one, from $\Sigma_n$ to $\Sigma_{n-1}$.
\end{dfn}

In~\cite{ea:mm}, the first author computed the fundamental group of Flenner-Za{\u\i}denberg curves and 
showed when it is non-abelian using orbifold groups. We show here that this can also be geometrically proved. 
In order to construct these curves, we start with a smooth conic $C$ with two tangent lines~$L_1$ and $L_2$,
intersecting at some point~$P$. After blowing-up~$P$ one obtains $\pi:\Sigma_1\to\bp^1$ with exceptional
component~$E$. Let $L_3$ be another line in the pencil through~$P$ which intersects~$C$ at two points~$Q_1$ 
and $Q_2$. Let us fix two positive integers $a,b$. After performing $a$ negative elementary Nagata 
transformations at the point corresponding to the fiber of~$L_1$ and $b$ at the point corresponding 
to the fiber of~$L_2$ one obtains a Hirzebruch surface $\Sigma_{a+b+1}$. One can then perform $a+b$ 
negative elementary Nagata transformations at the point~$Q_2$ (it corresponds to the fiber~$L_3$). 
After this process, $E$ can be blown down which turns our surface into~$\bp^2$. The curve~$C_{a,b}$ 
obtained has degree~$d:=a+b+2$ and three singular points of type~$\ba_{2a}$, $\ba_{2b}$, and a third one 
with local equation~$u^{d-2}=v^{d-1}$.

\begin{prop}
The fundamental group of $\bp^2\setminus C_{a,b}$ possesses a geometric surjection onto $\bp^1_{2,a+b,c}$, 
where $c:=\gcd(2 a+1,2 b+1)$.
\end{prop}

\begin{proof}
It is enough to follow the pencil of conics generated by $L_1+L_2$ and $C$ through the above transformations. 
We obtain a pencil of curves of degree~$2(d-1)$, where one fiber is $(2 a+1) \tilde{L}_1+(2 b+1) \tilde{L}_2$
(they are the lines corresponding to the fibers of $L_1$ and $L_2$. The fiber containing $C_{a,b}$ is of the form
$C_{a,b}+(d-2)\tilde{L}_3$. Finally the double line in the pencil becomes a double curve of degree~$d-1$.
\end{proof}

In~\cite{tono:01}, K.~Tono describes all rational unicuspidal curves such that its complement in~$\bp^2$ has
logarithmic Kodaira dimension~$1$. The construction given in~\cite[Theorem~1]{tono:01} shows that the complement 
of these curves have non-abelian fundamental group. Any other known rational unicuspidal curve has abelian
fundamental group (for the complement).

\begin{prop}
For any Tono's curve~$C$ their fundamental group possesses a geometric surjection onto $\bp^1_{\mu_A,\mu_G,n(C)}$,
where $\mu_A,\mu_G\geq 2$ and the number $n(C)$ is the opposite of the self-intersection of the strict transform of
$C$ after the minimal embedded resolution of its unique singular point. This number is at least~$2$.
\end{prop}

\begin{proof}
It is enough to consider the construction of~\cite[Theorem~1]{tono:01} where a pencil is obtained with two multiple
fibers~$\mu_A A$ and~$\mu_G G$ and a reducible fiber of the form $C+n(C) B$, where $B$ is either a line (type~I) or a smooth conic (type~II).
\end{proof}

\begin{ejm}
The curves of type~I are parametrized by two integers $n,s\geq 2$. The curve~$C$ has degree $(n+1)^2 (s-1)+1$, where $n(C)=n$, $\mu_A=n+1$ and $\mu_G=(n+1)(s-1)+1$. For $n=s=2$, we obtain the multiplicities $2,3,4$; in fact, one can compute that this group is finite.
\end{ejm}

\section{Degtyarev curves}\label{sec-dgt}

Let us consider a projective Degtyarev curve, i.e., a plane projective curve of degree~$5$ 
such that $\sing(C)$ consists of three points, and for each point $P\in\sing(C)$ the germ 
$(C,P)$ is topologically equivalent to an $\ba_4$-singularity, i.e. with local equation 
$v^2-u^5=0$; note that in this case, the germs are also analytically equivalent.

Most of the following properties appear in~\cite{deg:90} and~\cite{nmb:86}, but we include for 
the sake of completeness.

\begin{prps}\label{prps-dgt}
Let $C\subset\bp^2$ be a projective Degtyarev curve. Then:
\begin{enumerate}
\enet{\rm(D\arabic{enumi})}
\item\label{prps-dgt1} The curve $C$ is irreducible.
\item\label{prps-dgt2} 
The tangent line $L$ of $C$ at a singular point $P$ satisfies $(L\cdot C)_P=4$.
\item\label{prps-dgt3} 
Two Degtyarev projective curves are projectively equivalent.
\item\label{prps-dgt4} 
The subgroup of projective transformations preserving $C$ is cyclic of order~$3$.
\item\label{prps-dgt5} 
The curve $C$ is autodual.
\end{enumerate}
\end{prps}

\begin{proof}
Since the three singular points are locally irreducible, \ref{prps-dgt1} is true.
For \ref{prps-dgt2}, note that $4\leq (L\cdot C)_P\leq 5$. Let us assume that $(L\cdot C)_P=5$;
considering $L$ as the line at infinity, $C\setminus L$ is an affine curve homeomorphic to $\bc$.
This case is discarded using Za{\u\i}denberg-Lin Theorem~\cite{zl} and \ref{prps-dgt2} results.

In order to prove~\ref{prps-dgt3}, there are two approaches. The direct approach consists of 
computing the equations of the curve~$C$ fixing the position of the singular points and some of 
their tangent lines. The second method is quite simple and worth describing here: Let $C_1,C_2$ 
be two projective Degtyarev curves. By B\'ezout's Theorem, the singular points are not aligned; 
and hence, after a projective transformation, one may assume that $\sing(C_1)=\sing(C_2)=:S$. 
Assuming that $S:=\{[1:0:0],[0:1:0],[0:0:1]\}$, one can perform a standard Cremona 
transformation~$\psi:\bp^2\dashrightarrow\bp^2$ based on the three singular points and defined
by $\psi([x:y:z])=[y z: x z: x y]$.

Geometrically, this rational map is obtained by blowing-up the three vertices of $S$ (obtaining a 
rational surface $X_\ell$) and then blowing down the strict transforms of the lines joining the 
points of $S$ (which have self-intersection $-1$ in $X_\ell$). One can easily compute that 
$\tilde{C_i}:=\psi(C_i)$ is a tricuspidal quartic. It is well known that there is only one tricuspidal 
quartic, up to projective transformation, therefore, after a suitable change of coordinates, one may 
assume $\tilde{C_1}=\tilde{C_2}=:\tilde{C}$, where $\tilde{C}$ is the curve with equation given
in~\eqref{eq-tricusp}. The tricuspidal quartic satisfies the following properties. Let
$\sing(\tilde{C})=\{P_1,P_2,P_3\}$; there are three points 
$Q_1^\ell,Q_2^\ell,Q_3^\ell\in\tilde{C}$, $\ell=1,2$ such that $P_i,Q_j^\ell,Q_k^\ell$ are aligned 
for all the possibilities with $\#\{i,j,k\}=3$. Let $\ca_\ell$ be the arrangements of curves given 
by $\tilde{C}$ and the lines joining $Q_i^\ell$ and $Q_j^\ell$.

The curve $\tilde{C}$ is parametrized as in~\eqref{par-tricusp} and the singular points 
$P_1=[0:1:0]$, $P_2=[1:0:0]$, and $P_3=[0:0:1]$ correspond to $[t:s]=[0:1], [1:1],$ and $[1:0]$. 
It is not hard to check that $A_\ell:=(\alpha_\ell,2+\alpha_\ell,-\alpha_\ell)$ are affine parameters 
of $(Q_1^\ell,Q_2^\ell,Q_3^\ell)$.
The last condition implies that $\alpha_\ell^2+\alpha_\ell-1=0$. If $\alpha_1=\alpha_2$ then $\ca_1=\ca_2$.

The group of projective transformations fixing $\tilde{C}$ is the group of the permutation of the coordinates.
The mapping $[x:y:z]\rightmap{\sigma}[x:z:y]$ induces $[t:s]\mapsto [s:t]$ in the parametrization, and $[x:y:z]\rightmap{\tau}[y:z:x]$ induces $[t:s]\mapsto [s:s-t]$.

Let us assume that $\alpha_1\neq\alpha_2$
Applying the projective transformation~$\sigma$, results into two operations on $A_\ell$: the permutation 
$(1,3)$, and the change of parameters. Thus, $\sigma(A_1)=(-\alpha_1^{-1},(\alpha_1+2)^{-1},\alpha_1^{-1})=A_2$,
which implies $\sigma(\ca_1)=\ca_2$.

Note that any projective transformation sending $\ca_1$ to $\ca_2$ lifts to an isomorphism $X_1\to X_2$ 
and this isomorphism induces a projective transformation of the source~$\bp^2$, hence~\ref{prps-dgt3} results.

In order to prove~\ref{prps-dgt4} one can use a similar argument on the projective transformations fixing $C$
(this last property was communicated to the authors by C.T.C.~Wall).

The property~\ref{prps-dgt5} follows from~\ref{prps-dgt2} and Pl\"ucker generalized formul{\ae}, 
see ~\cite{nmb:86}. More precisely, given a curve $D$ and a point $P\in D$, the order of the curve is the 
degree of its dual curve of~$\check{D}$: 
\[
\deg(\check{D})=\deg(D)(\deg(D)-1)-\sum_{P\in D}(\mu(C,P)-1+ m(C,P)).
\]
This formula implies that $\deg(\check{D})=5$. The dual of a singular point of type $\ba_4$ is either 
of the same type or of type $\be_8$ (in case the tangent line has multiplicity of intersection~$5$ with 
the curve at the singular point). Thus~\ref{prps-dgt5} holds.
\end{proof}

\begin{rem}
Note that any two projective Degtyarev curves are isotopic. 
Using the direct approach, we can give a symmetric equation:
\begin{equation*}
\begin{split}
\left( 7+3\,\sqrt{5}  \right) ({x}^{3} {z}^{2}+x^2 y^3+ y^2 z^3)+ 
\left( 2\,\sqrt{5} +6 \right) ({x}^{3} y z+x y^3 z+ x y z^3)&+\\
+2 ({x}^{3}{y}^{2}+{x}^{2}{z}^{3}+y^3 z^2)+ 
\left( 33+11\,\sqrt{5}  \right) ({x}^{2} y {z}^{2}+ {x}^{2}{y}^{2}z+x{y}^{2}{z}^{2})=&0. 
\end{split}
\end{equation*}
Note that the permutation of two variables comes from the Galois transformation in $\bq(\sqrt{5})$. 
The curve also admits an equation with rational coefficients; in that case one of the singular points 
has rational coordinates but the other two are conjugate in $\bq(\sqrt{5})$:
\begin{equation}\label{eq-rat}
{z}^{2}{y}^{3}-z ( 33 xz+2 {x}^{2}+8 {z}^{2} ) {y}^{2}+ 
( 21 {z}^{2}+21 xz-{x}^{2} )  ( {z}^{2}+11xz-{x
}^{2} ) y+ ( x-18 z)  ({z}^{2}+11 xz-{x}^{2
}) ^{2}=0 
\end{equation}
\end{rem}

Properties~\ref{prps-dgt} imply that the affine Degtyarev curve is also \emph{rigid},
i.e. any two affine Degtyarev curves are projectively equivalent, and in particular, they are isotopic. 
In order to study its complement, it is convenient to assume that the line corresponds to the line at 
infinity and hence it is enough to consider the complement of the affine curve whose equation is obtained from~\eqref{eq-rat} by taking~$z=1$.

\begin{figure}
\centering
\includegraphics[scale=.5]{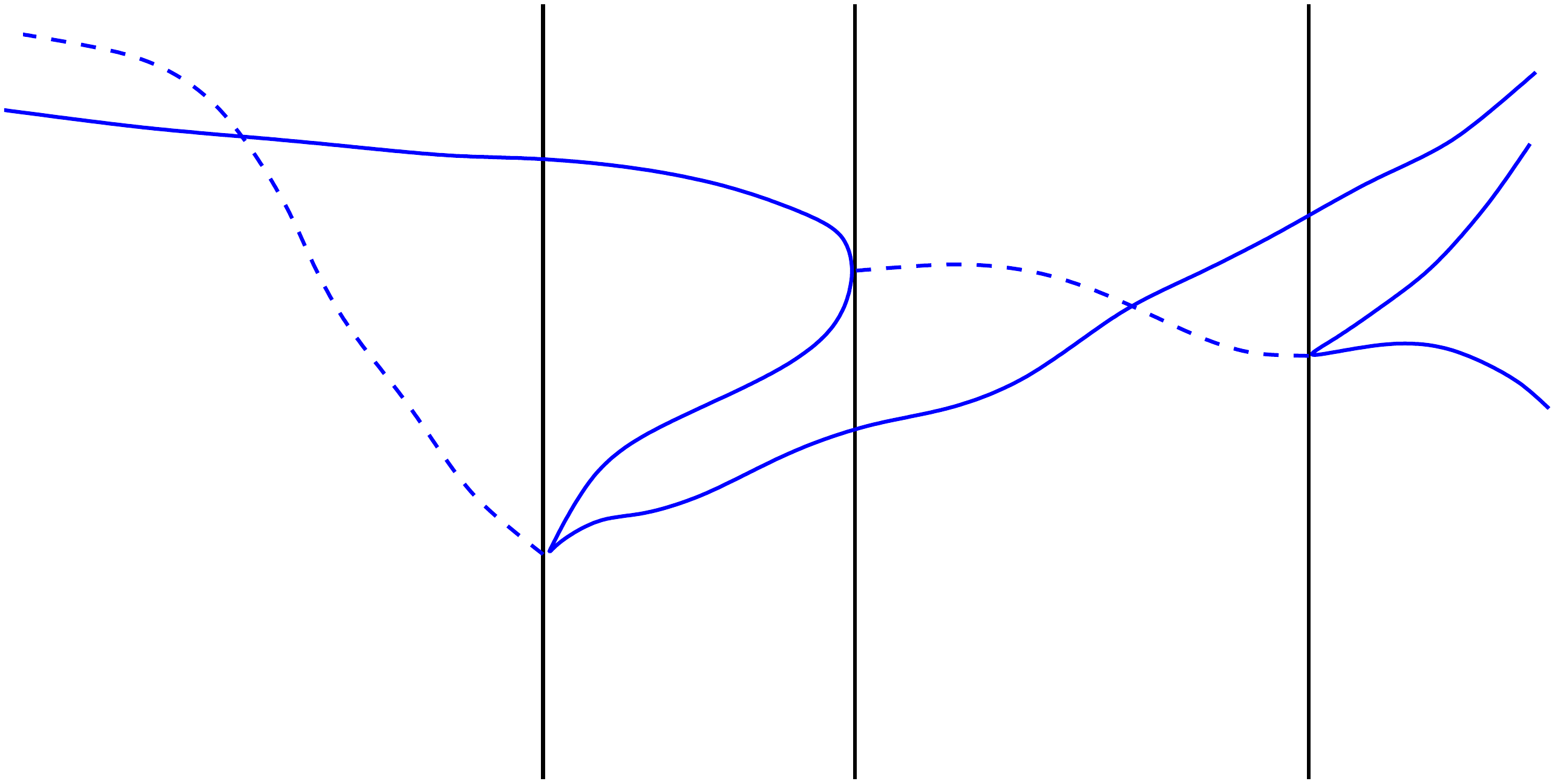}
\begin{picture}(0,0)
\put(-54,10){$x=\frac{11+5\sqrt{5}}{2}=a_+$}
\put(-54,180){$F_+$}
\put(-163,10){$x=0$}
\put(-163,180){$F_0$}
\put(-325,10){$x=\frac{11-5\sqrt{5}}{2}=a_-$}
\put(-260,180){$F_-$}
\end{picture}
\caption{Real picture of the affine Degtyarev curve}
\label{fig-a4a4a4}
\end{figure}
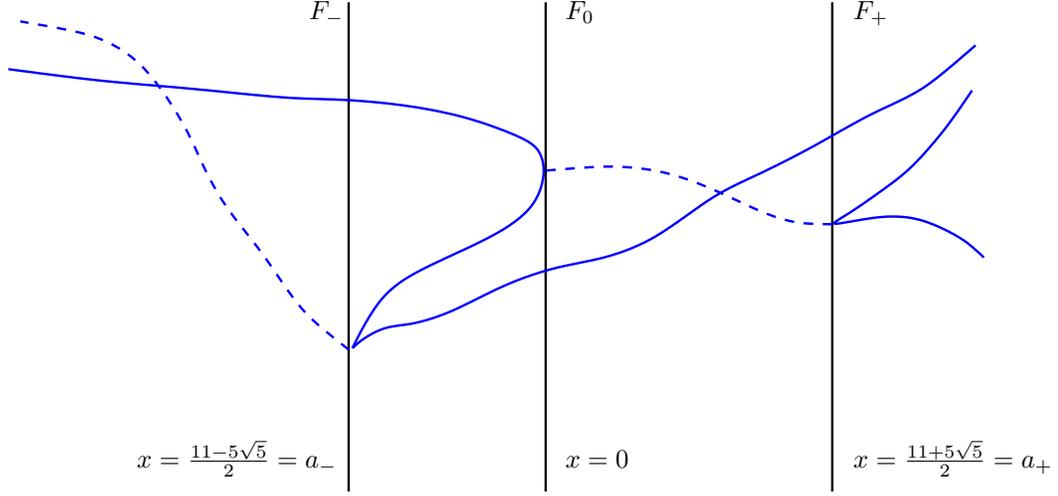

The fundamental group of the projective Degtyarev curve was computed in~\cite{deg:90}. 
Here we will compute the fundamental group of the affine curve and also show how to recover the group 
of the projective group. In order to compute the group we will use the braid monodromy associated with the 
projection $(x,y)\mapsto x$. Note that the discriminant of the equation~\eqref{eq-rat} (with $z=1$) is 
(up to a constant) $x (x^2-11 x-1)^5$. Since the three roots are real and the projection is $3:1$ 
with \emph{enough} real roots, the real picture in Figure~\ref{fig-a4a4a4} contains all the required
information to obtain the braid monodromy (the dotted lines represent the real part of the complex 
conjugate roots). 

The braid monodromy is defined as a representation $\nabla_0:\pi_1(\bc\setminus\{0,a_+,a_-\};*_0)\to\bb_3$. 
The source is a free group of rank three generated by:
\[
\mu_+:=\alpha_+\cdot\beta_+\cdot\gamma_+\cdot\alpha_+^{-1},\quad
\mu_0:=\alpha_+\cdot\beta_+\cdot\alpha_0\cdot\beta_0\cdot\gamma_0\cdot\alpha_0^{-1}\cdot\beta_+^{-1}\cdot\alpha_+^{-1}
\text{ and}
\]
\[
\mu_-:=\alpha_+\cdot\beta_+\cdot\alpha_0\cdot\beta_0\cdot\alpha_-
\cdot\beta_-\cdot\gamma_-\cdot
\alpha_-^{-1}\cdot\beta_0^{-1}\cdot\alpha_0^{-1}\cdot\beta_+^{-1}\cdot\alpha_+^{-1}.
\]
\begin{figure}
\centering
\includegraphics[scale=.7]{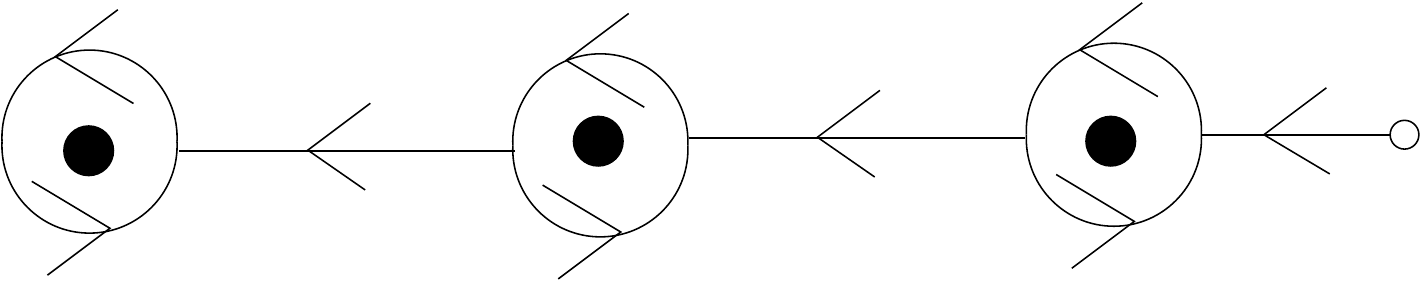}
\begin{picture}(0,0)
\put(-2,26){$*_0$}
\put(-30,10){$\alpha_+$}
\put(-60,50){$\beta_+$}
\put(-60,5){$\gamma_+$}
\put(-130,10){$\alpha_0$}
\put(-160,50){$\beta_0$}
\put(-160,5){$\gamma_0$}
\put(-230,10){$\alpha_-$}
\put(-260,50){$\beta_-$}
\put(-260,5){$\gamma_-$}
\end{picture}
\caption{Paths in $\bc\setminus\{0,a_+,a_-\}$}
\label{fig-base}
\end{figure}
Figure~\ref{fig-base} shows a geometric basis of $\pi_1(\bc\setminus\{0,a_+,a_-\};*_0)$. 
The braids are obtained by considering the way the roots with respect to $y$ move when the 
parameters move along $x$. We follow these conventions:
\begin{enumerate}
\enet{($\bb$\arabic{enumi})}
\item In order to draw the braids we consider the projection onto the real axis.
\item When two points have the same real part, we perturb the projection such that
positive imaginary parts go to the right and negative imaginary parts go to the right.
\item Roots will be numbered from right to left.
\item The above conventions give a canonical way to identify open braids with closed braids.
\end{enumerate}

Using the standard Artin generators of the braid groups, the braids obtained from following the paths 
in $\bc\setminus\{0,a_+,a_-\}$ shown in Figure~\ref{fig-base} are presented in Table~\ref{tab-braids}.
\begin{table}
\begin{center}
\begin{tabular}{|l|l|}\hline
\hfill\textbf{Paths}\hfill& \hfill\textbf{Braids}\hfill\\\hline
$\alpha_+$ & 1\\\hline
$\beta_+$& $\sigma_2^2$\\\hline
$\gamma_+$ & $\sigma_2^3$\\\hline
$\alpha_0$ & $\sigma_1^{-1}\sigma_2$\\\hline
$\beta_0$ & 1\\\hline
$\gamma_0$ & $\sigma_1$\\\hline
$\alpha_-$ & 1\\\hline
$\beta_-$ & $\sigma_2^2$\\\hline
$\gamma_-$ & $\sigma_2^3$\\\hline
\end{tabular}
\end{center}
\caption{Braids}
\label{tab-braids}
\end{table}

\begin{prop}\label{prop-braid}
The braid monodromy for the chosen projection of the affine Degtyarev curve is given by:
\[
\nabla_0(\mu_+)=\sigma_2^5,\quad
\nabla_0(\mu_0)=(\sigma_2^2\sigma_1^{-1}\sigma_2)*\sigma_1,\quad
\nabla_0(\mu_-)=(\sigma_2^2\sigma_1^{-1}\sigma_2\sigma_1)*\sigma_2^5=
\sigma_2^2*\sigma_1^5,
\]
where $a*b:=a b a^{-1}$.
\end{prop}

\section{Groups of Degtyarev curves}\label{sec-groups}

In order to compute the fundamental groups we apply the Zariski-van Kampen method. 
Let us consider the vertical line $F$ of equation~$x=*_0$. The set $F\setminus C$ is of the form $\{*_0\}\times\bc\setminus\{y_1,y_2,y_3\}$, where $y_1,y_2,y_3\in\br$, identified with the second 
factor. We choose a \emph{big} real $*$ number as base point. The free group $\pi_1(F\setminus C;*)$ 
has a free basis $g_1,g_2,g_3$ constructed as in Figure~\ref{fig-base}.
The natural action of $\bb_3$ on the free group $\bbf_3$ is expressed in this case as
\begin{equation}\label{action}
g_i^{\sigma_j}:=
\begin{cases}
g_{i+1}&\text{ if }i=j,\\
g_{i+1}*g_i&\text{ if }i=j+1,\\
g_i&\text{ if }i\neq j,j+1.
\end{cases}
\end{equation}

\begin{prop}\label{prop-pres}
The fundamental group of the affine Degtyarev curve has a presentation
\begin{equation}\label{pres-aff}
\left\langle
g_1,g_2,g_3\left|\ 
g_i^{\nabla_0(\mu_j)}=g_i, i=1,2,3, j=-,0,+
\right.
\right\rangle.
\end{equation}
In this presentation, the meridian of the line at infinity is $\left(g_3(g_2 g_1)^2\right)^{-1}$.
In particular, a presentation for the projective Degtyarev curve is
\begin{equation}
\left\langle
g_1,g_2,g_3\left|\ 
\eqref{pres-aff},\quad g_3=(g_2 g_1)^{-2}
\right.
\right\rangle.
\end{equation}
\end{prop}

\begin{proof}
The first presentation is a consequence of the Zariski-van Kampen method by means of the braid monodromy.
In order to prove the second one may consider a small deformation of the vertical line~$F$. It will intersect the curve at five points. Three of them are close to $(*_0,y_i)$, $i=1,2,3$, and the other two ones lie in the \emph{real} branches which go \emph{faster} to infinity. The boundary of a big disk in this line is the inverse
of a meridian of the line at infinity.
\end{proof}

\begin{rem}\label{rem-pres}
Proposition~\ref{prop-pres} provides right presentations of the group, but they may be quite cumbersome 
to work with by hand. Even if one wants to work with them with computer programs, like \texttt{GAP}\cite{GAP4}, 
the presentations could be intractable. There are several ways around this problem
\begin{enumerate}
\enet{\rm(P\arabic{enumi})}
\item\label{P1} 
The presentation \eqref{pres-aff} works if we replace the braid monodromy $\nabla_0$ for a conjugate.
For example, conjugating the braids in Proposition~\ref{prop-braid} by $\sigma_2^{2}$ produces simpler 
braids and hence a simpler presentation of the group.
\item\label{P2} 
Instead of finding a good braid to perform the conjugation in~\ref{P1} by inspection, one can try to 
interpret this conjugation in a geometric way. Changing the base point in $\bc\setminus\{0,a_+,a_-\}$ might 
produce simpler braids. For example choosing a real number $*_1\in (a_-,0)$ as a base point, one obtains the 
following as braid monodromy (for the new generators of the group):
\begin{equation}\label{newbraid}
\tilde{\mu}_+\mapsto(\sigma_2^{-1}\sigma_1)*\sigma_2^5,\quad
\tilde{\mu}_0\mapsto\sigma_1,\quad
\tilde{\mu}_-\mapsto\sigma_2^5.
\end{equation}
These braids have been obtained by conjugation of the ones in Proposition~\ref{prop-braid} by $\sigma_2^{2}\sigma_1\sigma_2^{-1}$.
\item\label{P3} 
If $g$ is a meridian of the line at infinity obtained using a braid monodromy~$\nabla_0$,
then, for a braid monodromy $(\nabla_0)^\tau:=\tau^{-1}\nabla_0\tau=(\tau^{-1})*\nabla_0$, 
a meridian of the line at infinity is~$g^\tau$.
\item\label{P4} 
There is another geometric way to reduce the presentation. Note that among the relations 
$(g_j)^{\sigma_2^5}=g_j$, $j=1,2,3$, one only needs to keep the relation given by $j=2$. First of all, 
the relation for $j=1$ is trivial; secondly $(g_3 g_2)^\tau=g_3 g_2$ and hence one of them is redundant. 
In the general case, this can be summarized as follows:
\begin{itemize}
\item Let us consider the action \eqref{action} (replacing $3$ by~$n$) of $\bb_n$ on the free group with 
basis $g_1,\dots,g_n$; let us consider a braid $\tau\in\bb_n$ which can be decomposed as 
$\tau=\tau_1\cdot\dots\cdot\tau_r$, where $\tau_j$ involves only a set of $n_j$ consecutive strings and $n=\sum_{j=1}^r n_j$. Then, among the relations $g_j^{\tau}=g_j$, we only need to consider 
$s:=\sum_{j=1}^r (n_j-1)=n-r$, disregarding one for each block of strings. Let $J_\tau$ be the chosen 
subset of indices.
\item If $\beta=(\tau)^\sigma$, and $\tau$ can be decomposed as above, then the set of relations 
$g_j^{\beta}=g_j$, $j=1,\dots,n$, is equivalent to $(g_j^\sigma)^\tau=g_j^\sigma$, $j\in J$.
\end{itemize}
For example, in our case the presentation~\eqref{pres-aff} can be reduced to have~$3$ relators.
\end{enumerate}
\end{rem}

\begin{prop}
The group~$G$ of the affine Degtyarev curve has a presentation:
\begin{equation}\label{pres-afin-def}
\left\langle
x, y\left|
x y x y x=y x y x y,
[x,y x y^{-1} x y x y^{-1} x y]=1
\right.
\right\rangle
\end{equation}
A presentation of the group~$G_\bp$ of the projective Degtyarev curve is obtained from \eqref{pres-afin-def}
by adding $x^5=1$. It turns out that $G_\bp$ is a group of order~$320$ with the following properties:
\begin{enumerate}
\enet{\rm($G_{\bp}$\arabic{enumi})}
\item $G_\bp/G_\bp'$ is cyclic of order~$5$.
\item The center $Z(G_\bp)$ the Klein group of order~$4$.
\item The group $G/Z(G_\bp)$ is a semidirect product of $(\bz/2)^4$ by $\bz^5$, where
the action of a generator of $\bz^5$ cyclically permutes a generator system $h_1,\dots,h_5$
of order~$2$ elements of $(\bz/2)^4$ satisfying $\sum h_i\equiv 0$.
\end{enumerate}
\end{prop}

\begin{proof}
The presentation of $G$ is obtained using the braid monodromy~\ref{newbraid} and Remark~\ref{rem-pres}\ref{P4},
where $x=g_1,g_2$ and $y=g_3$; note that $x$ and $y$ are conjugate. In order to obtain the presentation of $G_\bp$ 
the relation of the line at infinity needs to be added. This is a \emph{complicated} product of five conjugates
of~$x$. If one types this presentation in~\texttt{GAP}, the output is that $G_\bp$ has order~$320$ and that $x$ 
is an element of order~$5$. Also according to~\texttt{GAP}, the order of the quotient of $G$ obtained by adding 
the relation~$x^5=1$ is~$320$. These facts give the presentation of the statement. The properties of $G_\bp$ are
either trivial or easily computed using~\texttt{GAP}.
\end{proof}

\begin{prop}\label{prop-afin}
The group $G_\bp$ possesses no geometric surjections.
\end{prop}

\begin{proof}
The only properties we need for this are the sizes of the group and its abelianization. Let us assume that
$G_\bp$ possesses a geometric surjection. Since it is finite, the orbifold should be spherical, i.e.
either $\bp^1_{2,2,n}$, $n\geq 3$, or $\bp^1_{2,3,m}$, $m=3,4,5$. Since the order of the orbifold group must 
divide $320$, the only possibilities are $(2,2,n)$, where $n|160$. The group is dihedral and its 
abelianization is either $\bz/2$ or $(\bz/2)^2$; since the abelianization of $G_\bp$ is of order~$5$ the 
result follows.
\end{proof}

We finish this section with the main result of this paper. We are going to compute the characteristic 
varieties of the complement of the affine Degtyarev curve and we will prove that these components cannot 
come from the characteristic varieties of an orbifold.

\begin{thm}\label{thm-main}
Let $\bt_G=\bc^*$ be the character torus of $G$. Then $\cv_1(G)$ is the set containing~$1$ and the $10$-th 
primitive roots of unity, whereas $\cv_2(G)=\emptyset$. Therefore there is no geometric surjection of $G$ 
onto an infinite orbifold group.
\end{thm}

Since finite group orbifolds do not have characteristic varieties, the following Corollary holds.
\begin{cor}
No irreducible component of $\cv_1(G)$ is obtained as the pull-back of an irreducible component of the $\cv_1(\Gamma)$ where $\Gamma$ is an orbifold group.
\end{cor}

\begin{proof}[Proof of Theorem~{\rm\ref{thm-main}}]
We are going to change the presentation~\eqref{pres-afin-def}, by taking a new generator $t$ satisfying $y=x t$:
\begin{equation}\label{pres-afin-ab}
\left\langle
x, t\left|\ 
x t x^2 t x= t x^2 t x^2 t,
[x,t x t^{-1} x t x t^{-1} x t]=1
\right.
\right\rangle
\end{equation}
It is clear that $1\in \cv_1(G)\setminus\cv_2(G)$ since the non-twisted homology has rank~$1$. Let us 
consider a non-trivial character $\xi\in\bt_G$, which is identified by the image $1\neq\zeta$ of a 
positive generator of $\bz$. One can associate a $CW$-complex with the presentation~\eqref{pres-afin-ab}
with one $0$-cell~$p$, two $1$-cells $x,t$ and two $2$-cells $A,B$ (corresponding to the relations).
Then, the complex~$C_*(X;\bc)^\xi$ with which to compute the twisted homology is
$$
0\longrightarrow\bc^2\overset{\partial_2}{\longrightarrow}\bc^2\overset{\partial_1}{\longrightarrow}\bc\longrightarrow 0.
$$
The matrix for $\partial_1$ is
$
\left(
\begin{smallmatrix}
\zeta-1\\
0
\end{smallmatrix}
\right)$.
In particular, $\dim\ker\partial_1=1$ and hence $\cv_2(G)=0$. The matrix for $\partial_2$ equals
$$
\begin{pmatrix}
0&0\\
1-\zeta+\zeta^2-\zeta^3+\zeta^4&(1-\zeta+\zeta^2-\zeta^3+\zeta^4)(\zeta-1)
\end{pmatrix}.
$$
The homology is non trivial if and only if the matrix vanishes and hence $\cv_1(G)$ is as in the statement.

Since we are working with the complement of an affine (hence projective) curve, if $G$ admits a geometric 
surjection onto an infinite orbifold group, the orbifold must be over a rational curve. Since the 
abelianization has rank~$1$, the rational curve must be either $\bc$ or $\bp^1$. Any dominant morphism 
with target $\bc$ can be considered as dominant on $\bp^1$ and we treat only this case.

One needs to consider only orbifolds over~$\bp^1$ whose fundamental groups are infinite, have cyclic 
abelianizations and admit the $10$-th primitive roots of unity in their characteristic varieties. In 
particular, the abelianization must be of the type~$\bz/n$, where~$10$ divides~$n$.

Any such orbifolds admit dominant morphisms onto $\bp^1_{2,5,10}$ and $\bp^1_{2,2,5,5}$.
The properties of $\cv_2$ allow us to discard~$\bp^1_{2,2,5,5}$, see Proposition~\ref{prop-2255}.

Let us assume that there is a geometric surjection onto the orbifold $\bp^1_{2,5,10}$. 
Proposition~\ref{prop-2510} does not provide a direct obstruction in terms of $\cv_1$. 
We know that the kernel of the abelianization map is the fundamental group $K_2$ of a compact Riemann 
surface of genus~$2$, see Proposition~\ref{prop-2510}.

Note that $(x y)^5=(x^2 t)^5$ is a central element and the group~$K$ generated by its element defines 
an injection in $G/G'$. Following~\cite{deg:2010}, if $G_0:=G/K$, the groups $G_0'$ and $G'$ are 
isomorphic and hence $G'$ is finitely presented. Using Reidemeister-Schreier method, we find the 
following presentation:
\begin{equation}\label{pres-afin-sub}
G'=\left\langle
t_0,t_1,t_2,t_3,t_4\left|\ 
t_{n+1} t_{n+3}= t_n  t_{n+2} t_{n+4},
B_n=B_{n+1}
\right.
\right\rangle,
\end{equation}
where $B_n:=t_n t_{n+1}^{-1}t_{n+2} t_{n+3}^{-1} t_{n+4}$ and $x*t_n=t_{n+1}$. 
Note that $x^{10}*t_n=t_{n+10}=A*t_n$, where $A:=t_n  t_{n+2} t_{n+4}  t_{n+6} t_{n+8}$ for any $n$. 
This guarantees that the above presentation is finite. Summarizing, one can deduce that the kernel 
$K_1$ of the epimorphism onto $\bz/10$ equals $\bz\times G'$. Note that the rank of $K_1$ equals~$5$ 
and the rank of $K_2$ equals~$4$, so no contradiction arises.

According to \texttt{GAP} the following quotients of the lower central series have ranks~$5$ and $16$ 
for $K_2$, and $2$ and $0$ (order $5$) for $K_1$ and hence such an epimorphism cannot exist.
\end{proof}

\section{Further properties of the affine Degtyarev curve}\label{sec-ell}

The affine Degtyarev curve is related with elliptic fibrations as follows.
In order to work in a projective setting, one can first consider the projective Degtyarev curve, and
fix a singular point~$P$. We will denote by~$L$ the tangent line of~$C$ at~$P$, and the remaining 
singular points by $P_{\pm}$. Let $\sigma:\Sigma_1\to\bp^2$ be the blow-up of~$P$ where $E$ denotes 
the exceptional component. Strict transforms will follow Convention~\ref{cvn-str}. 

Each generic fiber of $\Sigma$ intersects $C$ at three points. There are four exceptions; 
three of them can be seen in Figure~\ref{fig-a4a4a4} and they are denoted by $F_+$, $F_0$, and $F_-$. 
The fourth one is $L$, which intersects $C$ at two points: one is smooth and transversal and the other 
one is the infinitely near point of $P$ in $E$, which is of type~$\ba_2$. In order to separate $C$ and 
$E$ we perform a positive elementary Nagata transformation $\rho:\Sigma_1\dashrightarrow\Sigma_2$ on the 
fiber corresponding to $L$. The fiber which replaces~$L$ is denoted by $F_\infty$. Note that $F_\infty$ 
intersects~$C$ at two points: one of them corresponds to the blow-down of $L$ and the other one is a 
point with a generic tangency. In particular, the combinatorics of the 
intersections at $F_0$ and $F_\infty$ coincides.

\begin{rem}
Properties~\ref{prps-dgt} imply the rigidity of this arrangement of curves in~$\Sigma_2$. In particular, 
once the four fibers are ordered the cross-ratio of their images in $\bp^1$ provides an invariant of the
arrangement. The existence of an automorphism of $\Sigma_2$ preserving~$C$ and exchanging the two fibers
containing the singular points can be easily checked. As a consequence of the cross-ratio argument, the 
two tangent fibers must also be exchanged. This automorphism defines a birational map of $\bp^2$ which 
is related to the two solutions in $\bq(\sqrt{5})$ exhibited in the proof of 
Property~\ref{prps-dgt}\ref{prps-dgt3}.
\end{rem}

Let us consider the minimal resolution~$Z$ of the double covering of $\Sigma_2$ ramified at $C+E$. 
The ruling of $\Sigma_2$ induces a morphism $\rho:Z\to\bp^1$ such that the generic fiber is elliptic. 
The only singular fibers are the preimages of $F_+$, $F_-$ (of type~$I_5$ in Kodaira notation), $F_0$, 
and $F_\infty$ (of type~$I_1$). These elliptic fibrations have been extensively studied in~\cite{mp:86}. 
Once a section is fixed (e.g. the preimage of~$E$), the set of sections has an abelian group structure 
(inherited by the structure on the fibers) which is called the Mordell-Weil group. Note that the involution associated with the double covering is defined by taking the opposite. It is known that the Mordell-Weil 
group of $Z$ is cyclic of order~$5$. 

Let us consider a conic $C_1$ tangent to $C$ both at $P$ and at another singular point and transversal to the 
third singular point. The preimage of $C_1$ by the double covering has two irreducible components which are 
denoted by $E_1$ and $-E_1$: they are opposite sections in the Mordell-Weil group. Interchanging the two 
singular points, one obtains the remaining two sections $E_2$ and $-E_2$ of $Z$.

Let us recall that $G$ denotes the fundamental group of the complement of the affine Degtyarev curve,
i.e. $\bp^2\setminus (C\cup L)=\Sigma_2\setminus(C\cup E\cup L_\infty)$. 

\begin{rem}
Despite Proposition~\ref{prop-afin}, note that its affine version, $G=\pi_1(\bp^2\setminus (C\cup L))$ 
\emph{does} posses a geometric surjection onto the orbifold over $\bp^1_{2,2,5}$, since $G$ admits an 
epimorphism onto the dihedral group of order~$10$, see for instance~\cite{act:05a}.
\end{rem}

In order to construct this morphism, we may use the ideas in~\cite{tokunaga4}. 
The mapping is obtained by a pencil of rational curves of degree~$10$, with the following non-reduced fibers:
\begin{itemize}
\item A smooth conic $C_2$ of multiplicity~$5$ such that $(C\cdot C_2)_{P_+}=2$,
$(C\cdot C_2)_{P_-}=4$ and $(C\cdot C_2)_{P}=4$.
\item A quintic $C_5$ of multiplicity~$2$ such that $(C\cdot C_5)_{P_+}=5$ ($P_+$ is a smooth point
of $C_5$), $(C\cdot C_2)_{P_-}=10$ ($P_-$ is a singular point of $C_5$ of type~$\ba_4$), and 
$(C\cdot C_2)_{P}=10$ ($P$ is a singular point of $C_5$ of type~$\bd_6$).
\item The curve $C+L+2D_2$ where~$D_2$ is a smooth conic such that 
$(C\cdot D_2)_{P_+}=0$, $(C\cdot D_2)_{P_-}=5$, and $(C\cdot D_2)_{P}=4$.
\end{itemize}

We finish this section by describing some properties of the group~$G$. For a point $Q\in C$,
the local fundamental group $\pi_1^{\text{\rm loc}}(C,Q)$ of $C$ at~$Q$ is $\pi_1(\bb_Q\setminus C)$,
where $\bb_Q$ is a Milnor ball. The inclusion $\bb_Q\setminus C\hookrightarrow\bc^2\setminus C$ induces 
a conjugacy class of subgroups (since the base point is not fixed) which will be called the 
\emph{image of the local fundamental group}.

\begin{prop}
Let $P_{\pm}$ be the two singular points of the affine Degtyarev curve.
\begin{enumerate}
\enet{\rm(\alph{enumi})}
\item The images of the local fundamental groups at $P_+$ and $P_-$ are the whole group~$G$.
\item The center of $G$ contains an abelian free subgroup of rank~$2$.
\end{enumerate}
\end{prop}

\begin{proof}
The property about the image of the local fundamental group at $P_-$ is obvious from the presentation~\eqref{pres-afin-def}. For~$P_+$ it can be deduced using~\texttt{GAP}. 
As a consequence we obtain two central elements (the images of the central elements of the local 
fundamental groups). The last property can be deduced by studying some quotients of subgroups of~$G$.
\end{proof}

\begin{thebibliography}{10}

\bibitem{ara:97}
D.~Arapura, \emph{Geometry of cohomology support loci for local systems. {I}},
  J. Algebraic Geom. \textbf{6} (1997), no.~3, 563--597.

\bibitem{ea:ta}
E.~Artal, \emph{A curve of degree five with non-abelian fundamental group},
  Topology Appl. \textbf{79} (1997), no.~1, 13--29.

\bibitem{ea:mm}
\bysame, \emph{Fundamental group of a class of rational cuspidal curves},
  Manuscripta Math. \textbf{93} (1997), no.~3, 273--281.

\bibitem{ACM-prep}
E.~Artal, J.I. Cogolludo, and D.~Matei, \emph{Orbifolds and characteristic
  varieties}, In preparation.

\bibitem{act:05a}
E.~Artal, J.I. Cogolludo, and H.~Tokunaga, \emph{Pencils and infinite dihedral
  covers of {$\mathbb P\sp 2$}}, Proc. Amer. Math. Soc. \textbf{136} (2008),
  no.~1, 21--29 (electronic).

\bibitem{Be}
A.~Beauville, \emph{Annulation du {$H\sp 1$} pour les fibr{\'e}s en droites
  plats}, Complex algebraic varieties ({B}ayreuth, 1990), Lecture Notes in
  Math., vol. 1507, Springer, Berlin, 1992, pp.~1--15.

\bibitem{Bu}
N.~Budur, \emph{Unitary local systems, multiplier ideals, and polynomial
  periodicity of {H}odge numbers}, Adv. Math. \textbf{221} (2009), no.~1,
  217--250.

\bibitem{deg:90}
A.I. Degtyar{\"e}v, \emph{Isotopic classification of complex plane projective
  curves of degree $5$}, Leningrad Math. J. \textbf{1} (1990), no.~4, 881--904.

\bibitem{deg:2010}
\bysame, \emph{Plane sextics via dessins d'enfants}, Geom. Topol. \textbf{14}
  (2010), no.~1, 393--433.

\bibitem{Delzant}
T.~Delzant, \emph{Trees, valuations and the {G}reen-{L}azarsfeld set}, Geom.
  Funct. Anal. \textbf{18} (2008), no.~4, 1236--1250.

\bibitem{dim:li}
A.~Dimca, \emph{Singularities and topology of hypersurfaces}, Universitext,
  Springer-Verlag, New York, 1992.

\bibitem{Di4}
\bysame, \emph{Characteristic varieties and constructible sheaves}, Atti Accad.
  Naz. Lincei Cl. Sci. Fis. Mat. Natur. Rend. Lincei (9) Mat. Appl. \textbf{18}
  (2007), no.~4, 365--389.

\bibitem{fz}
H.~Flenner and M.G. Za{\u\i}denberg, \emph{On a class of rational cuspidal
  plane curves}, Manuscripta Math. \textbf{89} (1996), no.~4, 439--459.

\bibitem{GAP4}
The GAP~Group, \emph{{GAP -- Groups, Algorithms, and Programming, Version
  4.4}}, 2004, available at \verb+(http://www.gap-system.org)+.

\bibitem{li:01}
A.~Libgober, \emph{Characteristic varieties of algebraic curves}, Applications
  of algebraic geometry to coding theory, physics and computation (Eilat,
  2001), Kluwer Acad. Publ., Dordrecht, 2001, pp.~215--254.

\bibitem{mp:86}
R.~Miranda and U.~Persson, \emph{On extremal rational elliptic surfaces}, Math.
  Z. \textbf{193} (1986), 537--558.

\bibitem{nmb:86}
M.~Namba, \emph{Geometry of projective algebraic curves}, Marcel Dekker Inc.,
  New York, 1984.

\bibitem{nem:87}
A.~N{\'e}methi, \emph{On the fundamental group of the complement of certain
  singular plane curves}, Math. Proc. Cambridge Philos. Soc. \textbf{102}
  (1987), no.~3, 453--457.

\bibitem{oka:75}
M.~Oka, \emph{Some plane curves whose complements have non-abelian fundamental
  groups}, Math. Ann. \textbf{218} (1975), no.~1, 55--65.

\bibitem{Si1}
C.~Simpson, \emph{Subspaces of moduli spaces of rank one local systems}, Ann.
  Sci. {\'E}cole Norm. Sup. (4) \textbf{26} (1993), no.~3, 361--401.

\bibitem{tokunaga4}
H.~Tokunaga, \emph{Dihedral coverings of algebraic surfaces and their
  application}, Trans. Amer. Math. Soc. \textbf{352} (2000), no.~9, 4007--4017.

\bibitem{tono:01}
K.~Tono, \emph{Rational unicuspidal plane curves with {$\overline\kappa=1$}},
  S\=urikaisekikenky\=usho K\=oky\=uroku (2001), no.~1233, 82--89, Newton
  polyhedra and singularities (Japanese) (Kyoto, 2001).

\bibitem{zl}
M.G. Za{\u\i}denberg and V.Ya. Lin, \emph{An irreducible, simply connected
  algebraic curve in {${\bf C}^{2}$} is equivalent to a quasihomogeneous
  curve}, Dokl. Akad. Nauk SSSR \textbf{271} (1983), no.~5, 1048--1052.

\bibitem{zr:29}
O.~Zariski, \emph{On the problem of existence of algebraic functions of two
  variables possessing a given branch curve}, Amer. J. Math. \textbf{51}
  (1929), 305--328.

\end{thebibliography}

\providecommand{\bysame}{\leavevmode\hbox to3em{\hrulefill}\thinspace}
\providecommand{\MR}{\relax\ifhmode\unskip\space\fi MR }
\providecommand{\MRhref}[2]{%
  \href{http://www.ams.org/mathscinet-getitem?mr=#1}{#2}
}
\providecommand{\href}[2]{#2}

\end{document}